# Groupoids: Unifying Internal and External Symmetry

### A Tour Through Some Examples

Alan Weinstein*
Department of Mathematics
University of California
Berkeley, CA 94720 USA
(alanw@math.berkeley.edu)

January 23, 1996

## 1  Introduction

Mathematicians tend to think of the notion of symmetry as being virtually synonymous with the theory of *groups* and their actions, perhaps largely because of the well known Erlanger program F. Klein and the related work of S. Lie, which virtually defined geometric structures by their groups of automorphisms. (See, for example, Yaglom's account in [21].) In fact, though groups are indeed sufficient to characterize homogeneous structures, there are plenty of objects which exhibit what we clearly recognize as symmetry, but which admit few or no nontrivial automorphisms. It turns out that the symmetry, and hence much of the structure, of such objects can be characterized algebraically if we use *groupoids* and not just groups.

*Research partially supported by NSF Grant DMS-93-09653. This article is based on a talk given first at UCB (the University of California at Berkeley) in April 1995, during Mathematics Awareness Week, the theme of which was Mathematics and Symmetry. I would like to thank the Department of Mathematics of another UCB (the University of Colorado at Boulder) for an invitation to Boulder in November 1995, where a similar talk was most recently given, and where the first draft of this manuscript was largely written in the downtown cafes. Thanks should also be given to a third UCB (Université Claude-Bernard, Lyon), where I started my intensive work on groupoids in the summer of 1986. Finally, I would like to thank Paul Brown for preparing the figures.

The aim of this paper is to explain, mostly through examples, what groupoids are and how they describe symmetry. We will begin with elementary examples, with discrete symmetry, and end with examples in the differentiable setting which involve *Lie groupoids* and their corresponding infinitesimal objects, *Lie algebroids*.

## 1.1 Some history

The following historical remarks are not intended to be complete but merely to indicate the breadth of areas where groupoids have been used. An extensive survey of groupoids as of 1986 can be found in R. Brown's article [2].

Groupoids were introduced (and named) by H. Brandt [1] in a 1926 paper on the composition of quadratic forms in four variables. C. Ehresmann [5] added further structures (topological and differentiable as well as algebraic) to groupoids, thereby introducing them as a tool in differentiable topology and geometry.

In algebraic geometry, A. Grothendieck used groupoids extensively and, in particular, showed how they could be used to tame the unruly equivalence relations which arise in the construction of moduli spaces. As the role of moduli spaces expands in physics as well as mathematics, groupoids continue to play an essential role.

In analysis, G. Mackey [11] used groupoids under the name of *virtual groups* to allow the treatment of ergodic actions of groups, "as if" they were transitive, while it was observed that the convolution operation, extended from groups to groupoids, made possible the construction of a multitude of interesting noncommutative algebras. (See [16] for a survey.) In this context, the use of groupoid convolution algebras as a substitute for the algebras of functions on badly behaved quotient spaces is a central theme in the noncommutative geometry of A. Connes [4]. The book just cited also shows the extent to which groupoids provide a framework for a unified study of operator algebras, foliations, and index theory.

In algebraic topology, the *fundamental groupoid* of a topological space has been exploited by P. Higgins, R. Brown, and others (see Brown's textbook [3]) where the use of a fixed base point as imposed by the usual fundamental group would be too restrictive. Groupoid methods are thus well adapted to disconnected spaces (which arise frequently when connected spaces are cut into pieces for study) and to spaces with fixed-point-free group

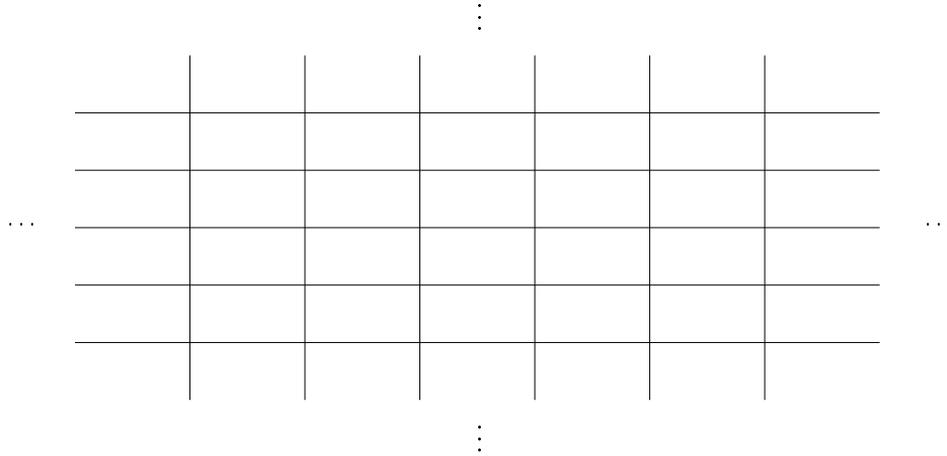

Figure 1: Tiling of the plane by $2 \times 1$ rectangles.

actions.

The extension of Lie theory from differentiable groups to groupoids was carried out for the most part by J. Pradines, as described in a series of notes ending with [15]. An exposition of this work, together with a detailed study of the role of groupoids and Lie algebroids in differential geometry, can be found in the book of K. Mackenzie [8].

Finally, the author's own interest in groupoids arose from the discovery [18][19], made independently at around the same time by Karasev [6] and Zakrzewski [22], that a groupoid with compatible *symplectic* structure is the appropriate geometric model for a family of noncommutative algebras obtained by deformation from the algebra of functions on a manifold with a Poisson bracket structure.

## 2 Global and local symmetry groupoids

We begin our exposition with a simple example which will lead to the definition of a groupoid. Consider a tiling of the euclidean plane $\mathbb{R}^2$ by identical $2 \times 1$ rectangles, specified by the set $X$ (idealized as 1-dimensional) where the grout between the tiles lies: $X = H \cup V$, where $H = \mathbb{R} \times \mathbb{Z}$ and $V = 2\mathbb{Z} \times \mathbb{R}$. (See Figure 1.) We will call each connected component of $\mathbb{R}^2 \setminus X$ a *tile*.

The symmetry of this tiling is traditionally described by the group $\Gamma$ consisting of those rigid motions of $\mathbb{R}^2$ which leave $X$ invariant. It consists of the normal subgroup of translations by elements of the lattice $\Lambda = H \cap V = 2\mathbb{Z} \times \mathbb{Z}$ (corresponding to corner points of the tiles), together with reflections through each of the points in $\frac{1}{2}\Lambda = \mathbb{Z} \times \frac{1}{2}\mathbb{Z}$ and across the horizontal and vertical lines through those points.

Now consider what is lost when we pass from $X$ to $\Gamma$.

- The same symmetry group would arise if we had replaced $X$ by the lattice $\Lambda$ of corner points, even though $\Lambda$ looks quite different from $X$. (Of course, there are reasons why this loss of detail in passing from $X$ to $\Gamma$ is welcome, as is the case for many mathematical abstractions, but read on.)

- The group $\Gamma$ retains no information about the local structure of the tiled plane, such as the fact that neighborhoods of all the points inside the tiles look identical if the tiles are uniform, while they may look different if the tiles are painted with a design (which could still be invariant under $\Gamma$).

- If, as is the case on real bathroom floors, the tiling is restricted to a finite part of the plane, such as $B = [0, 2m] \times [0, n]$, the symmetry group shrinks drastically. The subgroup of $\Gamma$ leaving $X \cap B$ invariant contains just 4 elements, even though a repetitive pattern on the bathroom floor is clearly visible.

Our first groupoid will enable us to describe the symmetry of the finite tiled rectangle. We first define the *transformation groupoid* of the action of $\Gamma$ on $\mathbb{R}^2$ to be the set[1]

$$G(\Gamma, \mathbb{R}^2) = \{(x, \gamma, y) | x \in \mathbb{R}^2, y \in \mathbb{R}^2, \gamma \in \Gamma, \text{ and } x = \gamma y\}$$

with the binary operation

$$(x, \gamma, y)(y, \nu, z) = (x, \gamma\nu, z).$$

This operation has the following properties.

1. It is defined only for certain pairs of elements: $gh$ is defined only when $\beta(g) = \alpha(h)$ for certain maps $\alpha$ and $\beta$ from $\Gamma$ onto $\mathbb{R}^2$; here $\alpha : (x, \gamma, y) \mapsto x$ and $\beta : (x, \gamma, y) \mapsto y$.

---

[1] As a set, $G(\Gamma, \mathbb{R}^2)$ is isomorphic to $\Gamma \times \mathbb{R}^2$, but we prefer the following more symmetric description, which makes the composition law more transparent.

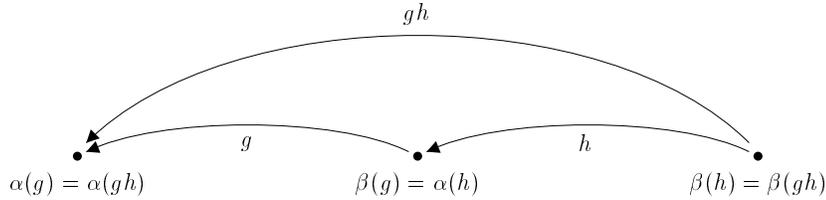

Figure 2: Multiplication in a groupoid.

2. It is associative: if either of the products $(gh)k$ or $g(hk)$ is defined, then so is the other, and they are equal.

3. For each $g$ in $G$, there are left and right identity elements $\lambda_g$ and $\rho_g$ such that $\lambda_g g = g = g\rho_g$.

4. Each $g$ in $G$ has an inverse $g^{-1}$ for which $gg^{-1} = \lambda_g$ and $g^{-1}g = \rho_g$.

More generally, a *groupoid* with base $B$ is a set $B$ with mappings $\alpha$ and $\beta$ from $G$ onto $B$ and a binary operation $(g, h) \mapsto gh$ satisfying the conditions 1-4 above. We may think of each element $g$ of $G$ as an arrow pointing from $\beta(g)$ to $\alpha(g)$ in $B$; arrows are multiplied by placing them head to tail, as in Figure 2.

Many properties of groupoids suggested by this picture are easily proven. for instance, $\alpha(g^{-1}) = \beta(g)$ (since $gg^{-1}$ is defined), and $\alpha(gh) = \alpha(g)$ (since $(\lambda_g g)h = \lambda_g(gh)$ implies that $\lambda_g(gh)$ is defined). One also shows easily that $\alpha(g)$ and $\lambda_g$ (or $\beta(g)$ and $\rho_g$) determine one another, thus producing a bijective mapping $\iota$ from the base $B$ to the subset $G^{(0)} \subseteq G$ consisting of the identity elements. The reader familiar with categories will probably have realized by now that *a groupoid is just a category in which every morphism has an inverse*.[2]

Now we return to our transformation groupoid $G(\Gamma, \mathbb{R}^2)$ and form its *restriction* to $B = [0, 2m] \times [0, n]$ (or any other subset $B$ of $\mathbb{R}^2$) by defining

$$G(\Gamma, \mathbb{R}^2)|_B = \{g \in G(\Gamma, \mathbb{R}^2) | \alpha(g) \text{ and } \beta(g) \text{ belong to } B\}.$$

The following concepts from groupoid theory, applied to this example, now show that $G(\Gamma, \mathbb{R}^2)|_B$ indeed captures the symmetry which we see in the tiled floor.

---

[2] Actually, not quite. The elements of a category can constitute a *class* rather than a set, so a groupoid is a *small* category with every morphism invertible.

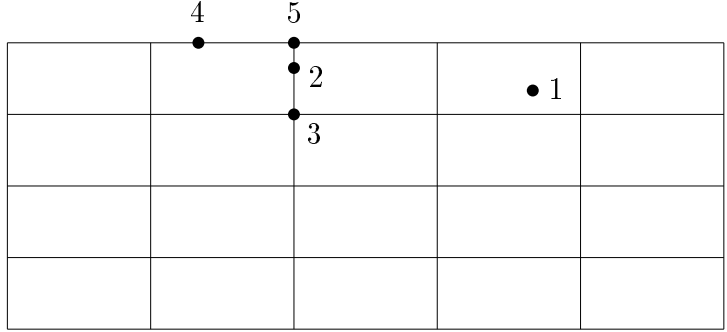

Figure 3: The point labeled by $j$ belongs to the orbit $\mathcal{O}_j$ of the groupoid $G(\Gamma,\mathbb{R}^2)|_B$.

- An *orbit* of the groupoid $G$ over $B$ is an equivalence class for the relation

  $x \sim_G y$ if and only if there is a groupoid element $g$ with $\alpha(g) = x$ and $\beta(g) = y$.

  In the example, two points are in the same orbit if they are similarly placed within their tiles or within the grout pattern. (We are allowed to turn the tiles over, as well as to translate them or rotate them by $180°$.)

- The *isotropy group* of $x \in B$ consists of those $g$ in $G$ with $\alpha(g) = x = \beta(g)$. In the example, the isotropy group is trivial for every point except those in $\frac{1}{2}\Lambda \cap B$, for which it is $\mathbb{Z}_2 \times \mathbb{Z}_2$.

We can describe local symmetries of our tiling by introducing another groupoid. Consider the plane $\mathbb{R}^2$ as the disjoint union of $P_1 = B \cap X$ (the grout), $P_2 = B \setminus P_1$ (the tiles), and $P_3 = \mathbb{R}^2 \setminus B$ (the exterior of the tiled room). Now define the local symmetry groupoid $G_{\mathrm{loc}}$ as the set of triples $(x,\gamma,y)$ in $B \times \Gamma \times B$ for which $x = \gamma y$, and for which $y$ has a neighborhood $\mathcal{U}$ in $\mathbb{R}^2$ such that $\gamma(\mathcal{U} \cap P_i) \subseteq P_i$ for $i = 1, 2, 3$. The composition is given by the same formula as for $G(\Gamma,\mathbb{R}^2)$.

For this groupoid, there are only a finite number of orbits (see Figure 3).

$\mathcal{O}_1$ = interior points of tiles
$\mathcal{O}_2$ = interior edge points

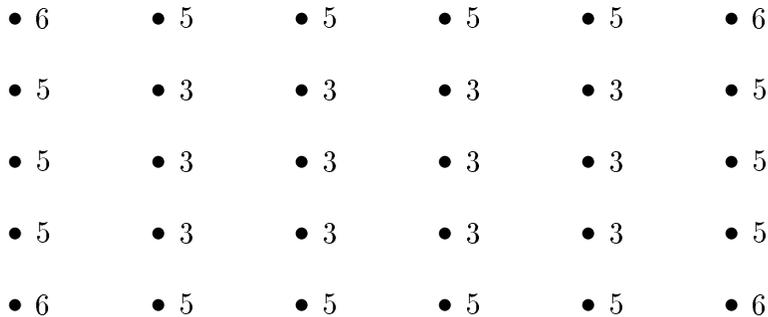

Figure 4: Orbit decomposition of the finite groupoid $G_{\text{loc}}|_L$.

$\mathcal{O}_3$ = interior crossing points
$\mathcal{O}_4$ = boundary edge points
$\mathcal{O}_5$ = boundary "T" points
$\mathcal{O}_6$ = boundary corner points

The isotropy group of a point in $\mathcal{O}_1$ is now isomorphic to the entire rotation group $O(2)$; it is $\mathbb{Z}_2 \times \mathbb{Z}_2$ for $\mathcal{O}_2$, the 8 element dihedral group $D_4$ for $\mathcal{O}_3$, and $\mathbb{Z}_2$ for $\mathcal{O}_4$, $\mathcal{O}_5$, and $\mathcal{O}_6$.

Now we can create a completely finite object by restricting $G_{\text{loc}}$ to the finite set $L = \mathcal{O}_3 \cup \mathcal{O}_5 \cup \mathcal{O}_6$ of tile corner points. $G_{\text{loc}}|_L$ now has just three orbits, with isotropy groups isomorphic to $D_4$, $\mathbb{Z}_2$, and $\mathbb{Z}_2$ as before. (See Figure 4.) At this point, we have something like the smile which remains after the Cheshire Cat has disappeared. The set $L$ has no structure at all except for that provided by the groupoid $G_{\text{loc}}|_L$, which is a sort of relic of the symmetry structure which was given by the original decomposition into grout, tiles, and and exterior. We could say that the points of $L$ have "internal symmetry" specified by the isotropy groups of $G_{\text{loc}}|_L$, and that the "types" of these points are classified by the orbit structure of the groupoid. This internal symmetry and type classification can be viewed either as a relic of the original tiling, or simply as a "geometric structure" on $L$.

The situation here is similar to that in gauge theory, where the points of space-time are equipped with internal symmetry groups whose representations correspond to physical particles. But gauge groups and gauge transformations are applicable only because all the points of space time look alike. To describe non-homogeneous structures like the tiled floor requires groupoids rather than groups.[3]

---

[3] Actually, groupoids are useful even in the realm of gauge theory. (See [9].)

A notable feature of the groupoid $G_{\text{loc}}|_L$ is that it can be entirely decomposed into its restrictions to the orbits $\mathcal{O}_3$, $\mathcal{O}_5$, and $\mathcal{O}_6$. In other words, the groupoid contains no information about the relation of the orbits to one another. This situation is typical for discrete groupoids. Only the presence of a topology on a groupoid and its base can interfere with such a decomposition.

## 3  Groupoids and equivalence relations

If $B$ is any set, the product $B \times B$ is a groupoid over $B$ with $\alpha(x,y) = x$, $\beta(x,y) = y$, and $(x,y)(y,z) = (x,z)$. The identities are the $(x,x)$, and $(x,y)^{-1} = (y,x)$. We call this the *pair groupoid* of $B$.[4] Note that a subgroupoid of $B \times B$, i.e. a subset closed under product and inversion, and containing all the identity elements,[5] is nothing but an equivalence relation on $B$. The orbits are the equivalence classes, and the isotropy subgroups are all trivial.

If $G$ is any groupoid over $B$, the map $(\alpha, \beta): G \to B \times B$ is a morphism from $G$ to the pair groupoid of $B$. (We leave it to the reader to formulate the precise definition of groupoid morphism.) The image of $(\alpha, \beta)$ is the orbit equivalence relation $\sim_G$, and the kernel is the union of the isotropy groups.

The morphism $(\alpha, \beta)$ suggests another way of looking at groupoids. Until now, we have been describing them as generalized groups, but now we want to think of groupoids as generalized equivalence relations as well. From this point of view, a groupoid over $B$ tells us not only *which* elements of $B$ are equivalent to one another (or "isomorphic"), but it also parametrizes the *different ways* ("isomorphisms") in which two elements can be equivalent. This leads us to the following guiding principle of Grothendieck, Mackey, Connes, Deligne, ....

> Almost every interesting equivalence relation on a space $B$ arises in a natural way as the orbit equivalence relation of some groupoid $G$ over $B$. Instead of dealing directly with the orbit space $B/G$ as an object in the category $\mathcal{S}_{\text{map}}$ of sets and mappings, one should consider instead the groupoid $G$ itself as an object in the category $\mathcal{G}_{\text{htp}}$ of groupoids and homotopy classes of morphisms.

---

[4] Other authors use the terms *coarse groupoid* or *banal groupoid* for this object.

[5] One often omits this last condition from the definition of a subgroupoid. To emphasize it, we can refer to a *wide subgroupoid*.

Here we will provide a definition. If

$$
\begin{array}{ccc}
G & \xrightarrow{f_G^i} & G' \\
\alpha \downarrow\downarrow \beta & & \alpha' \downarrow\downarrow \beta' \\
B & \xrightarrow{f_B^i} & B'
\end{array}
$$

are morphisms of groupoids for $i = 1, 2$, a *homotopy* (or *natural transformation* in the language of categories) from $f^1$ to $f^2$ is a map $h : B \to G'$ with the following properties.

- For each $x \in B$, $\alpha'(h(x)) = f_B^1(x)$ and $\beta'(h(x)) = f_B^2(x)$.

- For each $g \in G$, $h(\alpha(g))f_G^2(g) = f_G^1(g)h(\beta(g))$.

An isomorphism in $\mathcal{G}_{\text{htp}}$ between the groupoids $G$ and $G'$ (sometimes called an *equivalence of groupoids*) gives a bijection between the orbit spaces $B/G$ and $B/G'$, as well as group isomorphisms between corresponding isotropy subgroups in $G$ and $G'$. In the absence of further structure on the groupoids, an equivalence is essentially no more than this, but when the groupoids (and their bases) have measurable, topological, differentiable, algebraic, or symplectic (and Poisson) structures, the quotient spaces can be rather nasty objects, outside the category originally under consideration, and it is essential to focus attention on the groupoids themselves. One tool which makes this focus effective is an algebra associated to the groupoid $G$ over $B$ which is a useful substitute for the space of functions on $B/G$ when the quotient is a "bad" space. We describe this algebra in the next section.

## 4 Groupoid algebras

The *convolution* of two complex valued functions on a group $G$ is defined as the sum

$$(a * b)(g) = \sum_{k \in G} a(k)b(k^{-1}g),$$

at least when the support of each function (i.e. the set where it is not zero) is finite. When the supports are infinite, the sum may still be defined if the functions are absolutely summable. One obtains in this way, for instance, the rules for multiplying Fourier series as functions on the additive group $\mathbb{Z}$.

To extend the definition of convolution from groups to groupoids, one replaces the domain $G$ of the sum (or integral) by the $\alpha$-fibre $G_g = \alpha^{-1}(\alpha(g))$. Alternatively, one may use the more symmetric formula for convolution:

$$(a * b)(g) = \sum_{\{(k,\ell) | k\ell = g\}} a(k)b(\ell).$$

For instance, if $G = B \times B$ is the pair groupoid over $B = \{1, 2, \ldots, n\}$, the convolution operation (defined on all functions) is

$$(a * b)(i,j) = \sum_{(i,k)(k,j)=(i,j)} a(i,k)b(i,j),$$

and if we write $a_{ij}$ for $a(i,j)$ we get

$$(a * b)_{ij} = \sum_{k=1}^{n} a_{ik} b_{kj},$$

which exhibits $\mathcal{F}(\{1, 2, \ldots, n\} \times \{1, 2, \ldots, n\})$ as the algebra of $n \times n$ matrices.

This may look like a peculiar way of viewing the algebra of matrices, but Connes ([4] pp. 33-39, and elsewhere) has made the point that it was precisely as a groupoid algebra that Heisenberg constructed his original formulation of quantum mechanics. (The elements of the groupoid were transitions between excited states of an atom. Only later were the elements of the groupoid algebra identified as matrices.) The noncommutativity of the algebra of observables in quantum mechanics is then seen to be a direct consequence of the noncommutativity of the pair groupoid $B \times B$, as compared with the commutativity of the group $\mathbb{Z}$ underlying the Fourier series representation of oscillatory motion in classical mechanics.

When $G$ is a topological group(oid)[6] it is natural to look at continuous functions. Unless the $\alpha$-fibres of $G$ are discrete, though, there are no continuous functions which are summable along the fibres, and it necessary to replace the sums in the definition of convolution by integrals. For instance, if $G = \mathbb{R}^n \times \mathbb{R}^n$, we can define a convolution operation on the space $C_c(G)$ of continuous functions with *proper* support (i.e. the projection of the closure of $\{(x,y)|a(x,y) \neq 0\}$ onto each factor $\mathbb{R}^n$ is a proper map). Multiplication in $C_c(G)$ is given by the integral formula

$$(a * b)(x, y) = \int_{\mathbb{R}^n} a(x, z) b(z, y) \, dz, \tag{1}$$

where $dz$ is Lebesgue measure.[7]

---

[6]$G$ is a *topological groupoid* over $B$ if $G$ and $B$ are topological spaces, and $\alpha$, $\beta$, and multiplication are continuous maps.

This is precisely the composition law which we obtain by thinking of each $a \in C_c(G)$ as the (Schwartz) kernel of an operator $\tilde{a}$ on $L^2(\mathbb{R}^n)$:

$$(\tilde{a}\psi)(x) = \int_{\mathbb{R}^n} a(x, y)\psi(y) dy.$$

Not every operator (not even the identity) is realized by a continuous kernel, but all operators (including unbounded ones) satisfying reasonable continuity properties can be realized by kernels which are *distributions*, or generalized functions on $\mathbb{R}^n \times \mathbb{R}^n$.

To define convolution operations on functions on a general locally compact topological groupoid, we need to start with a family $\{\mu_x\}$ of measures along the $\alpha$-fibres, which will fill the role played by the measure $dz$ in (1) above. Convolution is then defined[8] by

$$(a * b)(g) = \int_{G_g} a(k) b(k^{-1} g) d\mu_g(k).$$

This operation is associative if the measures $\mu_g$ satisfy a suitable left-invariance property.

From continuous functions on $G$, on can pass by various processes of completion and localization to larger *groupoid algebras*, which play the role of continuous functions, measurable functions, or even distributions on the orbit space of $G$. We refer to the books [4], [14], and [16] for more details.

A more intrinsic construction of groupoid algebras, free of the choice of measures, can also be given. It includes as special cases the convolution of measures on a group $G$ and the multiplication of half-density kernels on $\mathbb{R}^n \times \mathbb{R}^n$. The definition can be found on p. 101 of [4] for the case of

---

[7] Dirac's notation $\langle x|a|y \rangle$ for $a(x, y)$ emphasizes the idea that the values of $a$ are "transition amplitudes" between "states" and so fits nicely with the groupoid interpretation of (1).

[8] This definition seems to be due to J. Westman [20].

Lie groupoids (see Section 5 below). A construction related to Mackey's intrinsic Hilbert spaces [10] should extend this idea to a general locally compact groupoid with a prescribed family of measure *classes* along the $\alpha$-fibres.

## 5 Lie groupoids and Lie algebroids

A groupoid $G$ over $B$ is a *Lie groupoid* if $G$ and $B$ are differentiable manifolds, and $\alpha$, $\beta$, and multiplication are differentiable maps.[9] For example, the pair groupoid $M \times M$ over any differentiable manifold $M$ is a Lie groupoid, and any Lie group is a Lie groupoid over a single point.

The infinitesimal object associated with a Lie groupoid is a *Lie algebroid*. A Lie algebroid over a manifold $B$ is defined to be a vector bundle $A$ over $B$ with a Lie algebra structure [ , ] on its space of smooth sections, together with a bundle map $\rho : A \to TB$ (called the *anchor* of the Lie algebroid) satisfying the conditions:

1. $[\rho(X), \rho(Y)] = \rho([X,Y])$;

2. $[X, \varphi Y] = \varphi[X,Y] + (\rho(X) \cdot \varphi)Y$.

Here, $X$ and $Y$ are smooth sections of $A$, $\varphi$ is a smooth function on $B$, and the bracket on the left hand side of the first condition is the usual Jacobi-Lie bracket of vector fields. For example, the tangent bundle $TB$ with the usual bracket, and the identity map as anchor, is a Lie algebroid over $B$, while any Lie algebra is a Lie algebra over a single point.

The Lie algebroid of a Lie groupoid is defined via left-invariant vector fields, just as it is in the case of Lie groups. Since the left translation by a groupoid element $g$ maps only from $G_{\beta(g)}$ to $G_{\alpha(g)}$, we must first restrict the values of our vector fields to lie in the subbundle $T^\alpha G = \ker(T\alpha)$ of $TG$. For every $x$ in $B$ and every $g$ in the $\beta$-fibre $\beta^{-1}(x)$, the fibre $T_g^\alpha G$ can be identified with $T_{\iota(x)}^\alpha G$ via left translation by $g$. Thus, the space of left-invariant sections of $T^\alpha G$ along $\beta^{-1}(x)$ can be identified with the single space $\mathcal{A}_x G = T_{\iota(x)}^\alpha G$, and $T\beta$ gives a well defined map $\rho_x : \mathcal{A}_x G \to T_x B$. The (disjoint) union of these objects for all $x$ in $B$ forms a vector bundle $\mathcal{A}G$ over $B$ with a bundle map $\rho : \mathcal{A}G \to TB$. (See Figure 5.) Identifying the sections of $\mathcal{A}G$ with left-invariant sections of $T^\alpha G \subseteq TG$ and using the usual bracket on sections of $TG$ defines a bracket of sections on $\mathcal{A}G$ making it into a Lie algebroid with anchor $\rho$.

---

[9]A further technical condition is included in the definition to insure that the domain of multiplication is a smooth submanifold of $G \times G$; the derivatives of $\alpha$ and $\beta$ are required to have maximal rank everywhere.

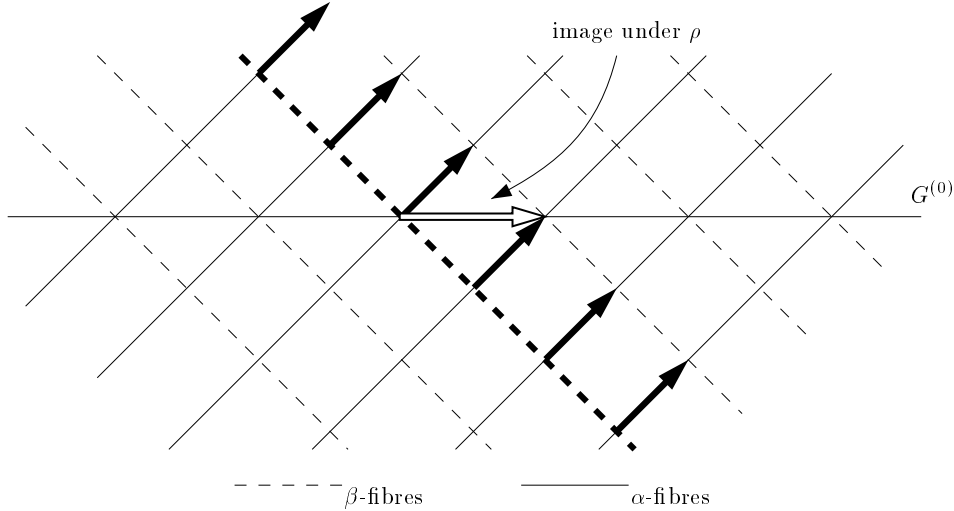

Figure 5: The groupoid $G$ with an element of $\mathcal{A}$ and its image under $\rho$.

For instance, the Lie algebroid of $B \times B$ is $TB$, and by extension of a well-known idea from Lie group theory, one may think of the Jacobi-Lie bracket on sections of $TB$ as being the infinitesimal remnant of the noncommutativity of the composition law $(x,y)(y,z) = (x,z)$.

The local versions of Lie's three fundamental theorems all extend from groups to groupoids, but there are some interesting differences at the global level. For instance, although every subalgebroid $A'$ of $\mathcal{A}G$ is indeed $\mathcal{A}G'$ for a Lie groupoid $G'$, there may be no such $G'$ for which the induced morphism $G' \to G$ is globally one-to-one. More seriously (and contrary to an incorrect assertion in [15]), there exist Lie algebroids which are "non-integrable" in the sense that they are not the Lie algebroid of any global Lie groupoid at all. (There is a *local* Lie groupoid corresponding to any Lie algebroid, though.) We refer to Chapter V of [8] for a discussion of all this, together with further references.

# 6  Boundary Lie algebroids

In this section, we will explain how Lie algebroids provide a natural setting for understanding some recent work by R. Melrose [12] [13] on the analysis of (pseudo)differential operators on manifolds with boundary. (Some prior knowledge about such operators on manifolds *without* boundary will be assumed here.)

The space $\mathcal{X}(M, \partial M)$ of vector fields which are tangent to the boundary $\partial M$ of a manifold $M$ forms a Lie algebra over $\mathbb{R}$ and a module over $C^\infty(M)$, and the condition $[X, \varphi Y] = \varphi[X, Y] + (X \cdot \varphi)Y$ is satisfied (since it is satisfied for all vector fields and functions). Remarkably, there is a vector bundle $^bTM$ over $M$ whose sections are in one-to-one correspondence with the elements of $\mathcal{X}(M, \partial M)$ via a bundle map $\rho : {}^bTM \to TM$ which is an isomorphism over the interior of $M$, and whose image over a boundary point is the tangent space to $\partial M$ at that point. These structures make $^bTM$ into a Lie algebroid. We refer to [12] for a precise definition of $^bTM$, mentioning here only that, when $(y_1, \ldots, y_{n-1}, x)$ are local coordinates for $M$ defined on an open subset of the half space $x \geq 0$, a local basis over $C^\infty(M)$ for the sections of $^bTM$ is given by $(\frac{\partial}{\partial y_1}, \ldots, \frac{\partial}{\partial y_{n-1}}, x\frac{\partial}{\partial x})$.

Melrose develops analysis on $M$ by regarding the $b$-tangent bundle $^bTM$ as the "correct" tangent bundle for this manifold with boundary. Thus, the algebra of $b$-differential operators is the algebra of operators on $C^\infty(M)$ generated by the sections of $^bTM$ (acting on functions via $\rho$) and $C^\infty(M)$ (acting on itself by multiplication).[10] The principal symbols of these operators are homogeneous functions on the dual, or $b$-cotangent bundle $^bT^*M$. Inverting (modulo smoothing operators) elliptic $b$-differential operators leads to the notion of $b$-pseudodifferential operators, whose symbols are again functions on the $b$-cotangent bundle.

Since the $b$-pseudodifferential operators act on $C^\infty(M)$, their Schwartz kernels are distributions on $M \times M$, which is a manifold with corners. Analysis of these kernels is facilitated, though, by lifting them to a certain manifold with boundary, $^b(M \times M)$, obtained from $M \times M$ by a blowing-up operation along the corners. The case where $M$ is a half-open interval is illustrated in Figure 6. The general case is locally the product of this one with a manifold without boundary.

It turns out that $^b(M \times M)$ is a groupoid whose Lie algebroid is $^bTM$, and that the composition law for $b$-kernels ([12], Equation (5.66)) is precisely the convolution operation in the groupoid algebra of ${}^b(M \times M)$. Thus, we can think of the groupoid ${}^b(M \times M)$ and its Lie algebroid ${}^bTM$ as providing a new "internal symmetry structure" appropriate to our manifold with boundary.

---

[10]Such a "universal enveloping algebra" was already defined for general Lie algebroids by Rinehart [17].

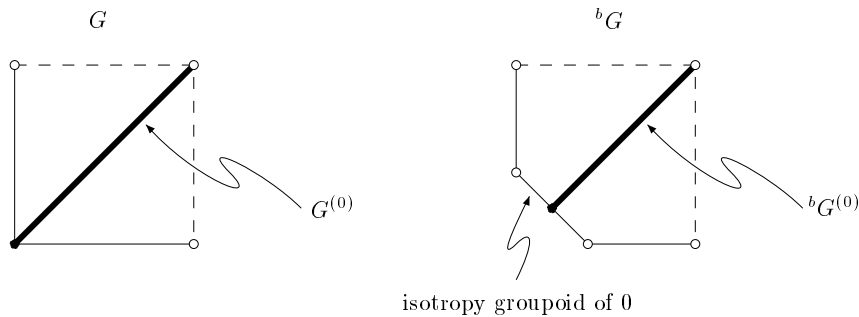

Figure 6: The groupoids $G = M \times M$ and ${}^bG = {}^b(M \times M)$, where $M = [0, 1)$.

It is therefore interesting to understand the structure of the groupoid ${}^b(M \times M)$. In fact, it decomposes algebraically (though not topologically) as the disjoint union of the pair groupoid over the interior of $M$ and the groupoid $G$ over $\partial M$ for which $(\alpha, \beta)^{-1}(x, y)$ consists of the orientation preserving linear isomorphisms between the normal lines $T_y M / T_y(\partial M)$ and $T_x M / T_x(\partial M)$. Thus, with the "geometric structure" on $M$ provided by the enlarged pair groupoid ${}^b(M \times M)$, each point of the boundary $\partial M$ has the internal symmetries are those of an oriented 1-dimensional vector space. We should think of the vectors in this space as being points "infinitesimally far" from the boundary, so that the boundary of $M$ becomes an "infinitesimal boundary layer." The section $x \frac{\partial}{\partial x}$ of ${}^bTM$ can now be seen as pointing into this layer in a non trivial way along $\partial M$; that is why it is not zero there.

At this point, our use of Lie algebroids and Lie groupoids adds only a geometric picture to Melrose's analysis. It does, however, suggest the possibility of extending the theory of pseudodifferential operators to algebras of operators associated with other Lie algebroids. In fact, the theory of longitudinal operators on foliated manifolds (see Section I.5.$\gamma$) of [4] is an example of such an extension, but the class of Lie algebroids which can be treated in this way remains to be determined. An interesting test case should be the Lie algebroids associated with the Bruhat-Poisson structures on flag manifolds [7], since the corresponding orbit decomposition on a flag manifold is the Bruhat decomposition, which plays a central role in representation theory.

## 7 Conclusion

This article has presented a small sample of the applications of groupoids in analysis, algebra, and topology. R. Brown has reported in [2] the suggestion of F.W. Lawvere that groupoids should perhaps be renamed "groups", and those special groupoids with just one base point given a new name to reflect their singular nature. Even if this is too far to go, I hope to have convinced the reader that groupoids are worth knowing about, and worth looking out for.[11]

## References

[1] Brandt, W., Über eine Verallgemeinerung des Gruppenbegriffes, *Math. Ann.* **96** (1926), 360-366.

[2] Brown, R., From groups to groupoids: a brief survey, *Bull. London Math. Soc.* **19** (1987), 113-134.

[3] Brown, R., *Topology: a geometric account of general topology, homotopy types, and the fundamental groupoid*, Halsted Press, New York, 1988.

[4] Connes, A., *Noncommutative Geometry*, Academic Press, San Diego, 1994.

[5] Ehresmann, C., Oeuvres complètes et commentées, ed. A.C. Ehresmann, Suppl. *Cahiers Top. Géom. Diff.*, Amiens, 1980-1984.

[6] Karasev, M.V., Analogues of objects of Lie group theory for nonlinear Poisson brackets, *Math. USSR Izvestia* **28**, (1987), 497-527.

[7] Lu, J.-H., and Weinstein, A., Poisson Lie groups, dressing transformations, and the Bruhat decomposition, *J. Diff. Geom.* **31** (1990), 501-526.

---

[11] Any reader stimulated to look further may wish to consult the *Groupoid Home Page*, http://amath-www.colorado.edu:80/math/department/groupoids/groupoids.shtml (under construction in January 1996).

[8] Mackenzie, K., *Lie Groupoids and Lie Algebroids in Differential Geometry*, LMS Lecture Notes Series, **124**, Cambridge Univ. Press, 1987.

[9] Mackenzie, K., Classification of principal bundles and Lie groupoids with prescribed gauge group bundle, *J. Pure Appl. Algebra* **58** (1989), 181-208.

[10] Mackey, G.W., *The Mathematical Foundations of Quantum Mechanics*, W.A. Benjamin, New York, 1963.

[11] Mackey, G.W., Ergodic theory and virtual groups, *Math. Ann.* **166** (1966), 187-207.

[12] Melrose, R.B., *The Atiyah-Patodi-Singer Index Theorem*, A.K. Peters, Wellesley, 1993.

[13] Melrose, R.B., *Geometric Scattering Theory*, Cambridge Univ. Press, Cambridge, 1995.

[14] Moore, C.C. and Schochet, C., *Global Analysis on Foliated Spaces*, MSRI Publications **9**, Springer-Verlag, New York (1988).

[15] Pradines, J., Troisième théorème de Lie sur les groupoïdes différentiables, *C. R. Acad. Sc. Paris,* **267** (1968), 21 - 23.

[16] Renault, J., A groupoid approach to $C^*$ algebras, *Lecture Notes in Math.* **793** (1980).

[17] Rinehart, G.S., Differential forms on general commutative algebras, *Trans. Amer. Math. Soc.* **108**, 195-222.

[18] Weinstein, A., Symplectic groupoids and Poisson manifolds, *Bull. Amer. Math. Soc.* **16**, (1987), 101-104.

[19] Weinstein, A., Noncommutative geometry and geometric quantization, *Symplectic Geometry and Mathematical Physics: actes du colloque en l'honneur de Jean-Marie Souriau*, P. Donato et al eds., Birkhäuser, Boston, 1991, pp. 446-461.

[20] Westman, J., Harmonic analysis on groupoids, *Pacific J. Math.* **27** (1968), 621-632.

[21] Yaglom, I. M., Felix Klein and Sophus Lie: evolution of the idea of symmetry in the Nineteenth Century, Birkhäuser, Boston, 1988.

[22] Zakrzewski, S., Quantum and classical pseudogroups, I and II, *Comm. Math. Phys.* **134** (1990), 347-370, 371-395.